\documentclass{article}
\usepackage{graphicx}
\RequirePackage[OT1]{fontenc}
\RequirePackage{amsthm,amsmath}
\RequirePackage[colorlinks]{hyperref}

\newcounter{counter_temp}

\theoremstyle{plain}
\newtheorem{thm}{Theorem}[section]
\newtheorem{lem}[thm]{Lemma}
\newtheorem{prop}[thm]{Proposition}
\newtheorem{ex}[thm]{Example}
\newtheorem{rem}[thm]{Remark}

\newcommand{\ve}[0]{\varepsilon}

\newcommand{\bbZ}{{\rm \rlap Z\kern 2.2pt Z}}

\newcommand{\Cov}{{\rm Cov}}
\newcommand{\Var}{{\rm Var}}
\newcommand{\Exp}{{\rm E}}

\newcommand{\MISE}{{\rm MISE}}
\newcommand{\CV}{{\rm CV}}

\newcommand{\ASE}{{\rm ASE}}

\newcommand{\convdistr}{\stackrel{\mbox{\rm d}}{\to}}
\newcommand{\convprob}{\stackrel{\mbox{\small\rm P}}{\to}}

\begin{document}

\title{Some results on random design regression with long memory errors and predictors\protect}

\author{Rafa{\l} Kulik \\
 Department of Mathematics \& Statistics,  University of Ottawa \\ 585 King Edward Avenue, Ottawa ON K1N 6N5, Canada \\ email: rkulik@uottawa.ca
\and \\
Pawe{\l} Lorek \\
Mathematical Institute, University of Wroc\l aw\\ pl. Grunwaldzki 2/4, 50-384 Wroc\l aw, Poland \\
email: lorek@math.uni.wroc.pl
}
\date{June 14, 2010}

\maketitle

\begin{abstract}
This paper studies nonparametric regression with long memory (LRD)
errors and predictors. First, we formulate general conditions which
guarantee the standard rate of convergence for a nonparametric
kernel estimator. Second, we calculate the Mean Integrated Squared
Error (MISE). In particular, we show that LRD of errors may
influence MISE. On the other hand, an estimator for a shape function
is typically not influenced by LRD in errors. Finally, we
investigate properties of a data-driven bandwidth choice. We show
that Averaged Squared Error (ASE) is a good approximation of
$\MISE$, however, this is not the case for a cross-validation
criterion.
\end{abstract}

\section{Introduction}
Consider the random design regression model,
\begin{equation*}
Y_i=m(X_i)+\varepsilon_i, \qquad i=1,\ldots,n ,
\end{equation*}
where $X_i$, $\varepsilon_i$, $i=1,\ldots,n$, are two mutually
independent sequences of random variables. We investigate the
problem of estimating function $m(\cdot)$. This problem is quite
well understood for weakly dependent data. Also, in the last two
decades, it has received a lot of attention in case of long range
dependence (LRD). In this situation most of results have been
obtained under quite specific assumptions on the errors and/or
predictors. In particular, it is typically assumed that both
sequences are infinite order moving averages, or they are defined as
(nonlinear) functionals of Gaussian sequences. However, in the
recent years different authors proposed many new (nonlinear) LRD
models. Although it is reasonable to assume some structure on the
observable predictors, particular assumptions on the errors are
almost impossible to verify. Therefore, one of our goals is to state
general conditions, which guarantee appropriate limit theorems.
 \\

In Section \ref{sec:clt} we discuss central limit theorem for the
Nadaraya-Watson estimator of $m(\cdot)$. As it is well known, for
LRD data we have a dichotomous behaviour: if a bandwidth $h$ is {\it
small}, then the rate of convergence is $\sqrt{nh}$, the same as in
i.i.d. case. Otherwise, if bandwidth is \textit{large}, long memory
contributes. We refer the reader to \cite{MielniczukWu2004} for the
most general result; see also \cite{ChengRobinson1994},
\cite{CsorgoMielniczuk1999}, \cite{CsorgoMielniczuk2000}. We state
general conditions, which guarantee $\sqrt{nh}$ rate of convergence.
These conditions can be easily verified for (subordinated) linear
LRD processes, FARIMA-GARCH models, stochastic volatility models
(including LARCH), as well as for antipersistent errors. In
particular, if $\varepsilon_i$, $i\ge 1$, is a linear process such
that $\Var\left(\sum_{i=1}^n\varepsilon_i\right)\sim Cn^{2-\alpha}$,
$\alpha\in (0,1)$, then $\sqrt{nh}$-consistency holds if
$hn^{1-\alpha}\to 0$. This agrees with previous results. (Here and
in the sequel, $C$ is a generic constant). On the other hand,
however, if the errors are described by a stochastic volatility,
then $\sqrt{nh}$-consistency {\it always} holds.\\

To verify that the condition $hn^{1-\alpha}\to 0$ holds, we need to
know the parameter $\alpha$. However, random variables
$\varepsilon_i$, $i=1,\ldots,n$, are not directly observable, and a
performance of various estimators of $\alpha$ is not clear.
Therefore, we will modify our estimation problem. In some cases, for
a purpose of an exploratory data analysis, it suffices to estimate a
{\it shape function}, $m^*(x)=m(x)-\int mf$. As indicated in
\cite{Efromovich1999}, \cite{Yang2001} and
\cite{KulikRaimondo2009b}, LRD of errors does not influence
estimation of $m^*(\cdot)$. This effect is proven here for the Nadaraya-Watson estimator. In fact, for
the linear LRD processes mentioned above, we obtain $\sqrt{nh}$-consistency if the predictors are independent and $h^5n^{1-\alpha}\to 0$, which is much weaker then the previous condition. However, this approach does not work in case of LRD predictors.  \\

Next, we investigate properties of the Mean Integrated Squared Error
($\MISE$). We assume that $\varepsilon_i$, $i\ge 1$, is the linear
process as mentioned above. We show that the optimal $\MISE$ has a
dichotomous behaviour: if the memory parameter $\alpha$ is greater
than $4/5$, then the rate of convergence is $n^{-4/5}$, the same as
in for i.i.d. errors. However, if $\alpha<4/5$, then the rate of
convergence is $n^{-\alpha}$. Interestingly, a possible LRD of
predictors does not influence the asymptotic behaviour of $\MISE$.
Similar results were obtained for density estimation, see
\cite{HallHart1990a}. On the other hand, in a fixed-design case, LRD
always influences the rates of convergence. For details we refer to
\cite{HallHart1990}.

To reduce the influence of LRD on $\MISE$, we may consider the shape
function. It is shown that for independent predictors, the Mean
Integrated Squared Error corresponding to $m^*$, has the same
asymptotic behaviour as in case of i.i.d. errors. In other words,
from expected-risk point of view, long memory does not influence
estimation of the shape function. We also note in passing that the
optimal bandwidth choice for the shape estimation
agrees with the optimal one for the estimation of the original function, as long as $\alpha>2/5$.\\

With help of formulas for $\MISE$, we obtain the optimal bandwidths.
As usual, they are not quite practical, since they involve unknown
parameters and a data-driven method has to be used. As argued in
\cite{HallLahiriPolzehl1995}, a plug-in method has some advantages
over cross-validation. However, let us note that the optimal
bandwidth is of the form $Cn^{-1/5}$ if $\alpha>2/5$, and
$Cn^{-(1-\alpha)/3}$, otherwise. Consequently, we have to know
$\alpha$ to be able to construct an appropriate plug-in bandwidth
selector. Therefore, we focus on cross-validation. Let us indicate
first that cross-validation is the valid procedure in a fixed-design
case. The procedure produces a bandwidth which is {\it close} to the
optimal one, and cross-validation criterion itself is a good
approximation to $\MISE$. The reader is referred to
\cite{HallLahiriPolzehl1995}. In the density estimation case,
however, cross-validation is a good approximation to  $\MISE$ {\it
if and only if} $\alpha>4/5$, see \cite{HallLahiriTruong1995} and
\cite{HallClaeskens2002} as well as discussion in Section
\ref{sec:bandwidth} for more details.

We will show that for random-design regression, the empirical
minimizer of the cross-validation criterion is a good approximation
to the optimal $h$,
however, the cross-validation itself is the valid procedure for $\alpha>4/5$ (which agrees with findings in \cite{HallLahiriTruong1995} and \cite{HallClaeskens2002}).  \\

The paper is organized as follows. In Section \ref{sec:assumptions} we collect assumptions and define estimators. Sections \ref{sec:clt}, \ref{sec:mse} and \ref{sec:bandwidth} contain results on central limit theorem, mean square error and bandwidth choice, respectively. In Section \ref{sec:numerical} we illustrate our findings by simulations. Next, in Section \ref{sec:examples} we present examples of time series, where it is possible to verify our conditions. Finally, the proofs are presented in the last section.

\section{Results}
\subsection{Assumptions and estimators}\label{sec:assumptions}
We will consider the following assumptions on the predictors $X_i$, $i\ge 1$:
\begin{itemize}
\item[{\rm (P1)}] $X_i$, $i\ge 1$, is a sequence of i.i.d. random variables. In this case, let ${\cal X}_i=\sigma(X_i,\ldots,X_1)$
\item[{\rm (P2)}] $X_i$, $i\ge 1$, is the infinite order moving average
$$
X_i=\sum_{k=0}^{\infty}a_k\zeta_{i-k},\qquad a_0=1,
$$
where $\zeta_i$, $-\infty<i<\infty$, is the sequence of centered,
i.i.d. Gaussian random variables and for $\alpha_X\in (0,1)$,
$a_k=A_0k^{-(\alpha_X+1)/2}$ for some $A_0$. Consequently, $X_i$ are
Gaussian and they are assumed to have unit variance. In this case we
denote ${\cal X}_i=\sigma(\zeta_i,\zeta_{i-1},\ldots)$. Furthermore,
note that $\Var\left(\sum_{i=1}^nX_i\right)\sim
A_1^2n^{2-\alpha_X}$, where $A_1$ is a finite constant.
\end{itemize}
We note that from the point of view of our results stated below, (P1) can be treated as the special case of (P2), by plugging-in $\alpha_X=1$. Thus, the results which are stated under (P2) assumption are valid also under (P1).\\

Under (P1), we do not need to assume a particular structure of errors $\varepsilon_i$, $i\ge 1$. The general assumption is
\begin{itemize}
\item[{\rm (E0)}] $\varepsilon_i=G(\eta_i,\eta_{i-1},\ldots)$, $i\ge 1$, where $\eta_i$, $-\infty<i<\infty$ is an i.i.d. sequence, independent of $X_i$, $i\ge 1$. We also assume that $\Exp[\varepsilon_1]=0$, $\Exp[\varepsilon_1^2]=1$. It assures that $\varepsilon_i$, $i\ge 1$, is a stationary ergodic sequence. Denote ${\cal H}_{i}=\sigma(\eta_i,\eta_{i-1},\ldots)$.
\end{itemize}
Under (P2) we will assume that
\begin{itemize}
\item[{\rm (E2)}] $\varepsilon_i$, $i\ge 1$, is the infinite order moving average
$$
\varepsilon_i=\sum_{k=0}^{\infty}c_k\eta_{i-k},\qquad c_0=1,
$$
where $\eta_i$, $-\infty<i<\infty$, is the sequence of centered,
i.i.d. random variables with a finite fourth moment,
$\Exp[\varepsilon_1^2]=1$, independent of $X_i$, $i\ge 1$. We also
assume that for $\alpha\in (0,1)$, $c_k\sim C_0k^{-(\alpha+1)/2}$,
as $k\to\infty$. Then
$\Var\left(\sum_{i=1}^n\varepsilon_i\right)\sim C_1^2n^{2-\alpha}$,
where $C_1$ is a constant.
\end{itemize}
We consider the
classical Nadaraya-Watson estimator
\begin{equation}\label{eq:kernel-estim}
\hat m(x)=\hat m_h(x)=\frac{1}{nh}\frac{1}{
\hat f_h(x)}\sum_{i=1}^nY_iK\left(\frac{x-X_i}{h}\right),
\end{equation}
where a nonnegative and bounded kernel $K$ fulfills the usual conditions:
$$\int K(u)du=1,\quad \int uK(u)du=0,\quad \kappa_2:=\int u^2K(u)\; du.$$
For a future use we denote $K_h(\cdot)=K(\cdot/h)$ and $\kappa_1:=\int K^2(u)du\not=0$.
The bandwidth $h=h_n$ fulfills the usual conditions $h\to 0$ and $nh\to \infty$. Furthermore, $\hat f_h(x)$ is the standard kernel estimator of the density $f$ of $X_1$.

Define further the {\it shape function}
$$ m^*:=m-\int mf
$$
and its estimator $\widehat{m^*}=\hat m-\widehat{\int mf}$, where
$$
\widehat{\int mf}:=\frac{1}{n}\sum_{i=1}^nY_i.
$$
Note that the latter is the unbiased estimator of $\int mf$, i.e.
$$
\frac{1}{n}\sum_{i=1}^n\Exp\left[Y_i\right]=\Exp\left[m(X_1)\right]=\int
mf.
$$
Finally, we will assume that $f$ and $m$ are twice differentiable
with bounded derivatives and that $f(x)>0$ for each $x$.
\subsection{Central limit theorems}\label{sec:clt}
Throughout this section we assume that (E0) holds. Let us formulate
the following conditions:
\setcounter{counter_temp}{\value{equation}}
\renewcommand{\theequation}{\Alph{equation}}
\setcounter{equation}{0}
\begin{equation}\label{eq:negligibility-condition}
\frac{\sqrt{nh}}{n}\sum_{i=1}^n (\Exp[\varepsilon_i|{\cal H}_{i-1}]-\varepsilon_i)=o_P(1).
\end{equation}
\begin{subequations}
\renewcommand{\theequation}{\theparentequation\arabic{equation}} 
\begin{equation}\label{eq:h-negligible-1}
\frac{\sqrt{nh}}{n}\sum_{i=1}^n\Exp[\varepsilon_i|{\cal H}_{i-1}]=o_P(1).
\end{equation}
\begin{equation}\label{eq:h-negligible}
\frac{\sqrt{nh}h^2}{n}\sum_{i=1}^n\Exp[\varepsilon_i|{\cal H}_{i-1}]=o_P(1).
\end{equation}
\end{subequations}
\begin{subequations}
\renewcommand{\theequation}{\theparentequation\arabic{equation}}
\begin{equation}\label{eq:h-negligible-2a}
h^5n^{1-\alpha_X}\to 0.
\end{equation}
\begin{equation}\label{eq:h-negligible-2}
hn^{1-\alpha_X}\to 0.
\end{equation}
\end{subequations}
\setcounter{equation}{\value{counter_temp}}
\renewcommand{\theequation}{\arabic{equation}}
Define
$$m_n(x):=\Exp[K((x-X_1)/h)m(X_1)]/\Exp[K((x-X_1)/h)], \quad m_n^*(x)=m_n(x)-\int mf.$$
Clearly, the bias is
\begin{equation}\label{eq:bias}
{\rm bias}(x)=m_n(x)-m(x)\sim h^2 \left(\frac{m''(x)}{2}+\frac{m'(x)f''(x)}{f(x)}\right)\int u^2K(u)\; du.
\end{equation}

\begin{prop}\label{thm:main}
Assume {\rm (P1)} and {\rm (E0)}. Under the conditions {\rm
(\ref{eq:negligibility-condition})} and {\rm
(\ref{eq:h-negligible-1})} we have
\begin{equation}\label{eq:result-1}
\sqrt{\frac{nh}{\kappa_1}}\sqrt{\hat f(x)}\left(\hat
m(x)-m_n(x)\right)\convdistr N(0,1).
\end{equation}
Assume {\rm (P2)} and {\rm (E2)}. Under the conditions {\rm
(\ref{eq:negligibility-condition})}, {\rm (\ref{eq:h-negligible-1})}
and {\rm (\ref{eq:h-negligible-2a})}, the asymptotics {\rm
(\ref{eq:result-1})} holds.
\end{prop}
\begin{rem}{\rm
Note that there is no symmetry between LRD assumptions on
$\varepsilon_i$ and $X_i$. For example, assume that $\varepsilon_i$,
and thus $\Exp[\varepsilon_i|{\cal H}_{i-1}]$ is the linear process
with $\Var\left(\sum_{i=1}^n\Exp[\varepsilon_i|{\cal
H}_{i-1}]\right)\sim C_1^2n^{2-\alpha}$, $\alpha\in (0,1)$; see
Example \ref{ex:linear}. Then {\rm (\ref{eq:h-negligible-1})} holds
if $hn^{1-\alpha}\to 0$, whereas the assumption for $\alpha_X$
requires $h^5n^{1-\alpha_X}\to 0$. }
\end{rem}
\begin{prop}\label{thm:main-1}
Assume {\rm (P2)}. Under the conditions {\rm
(\ref{eq:negligibility-condition})}, {\rm (\ref{eq:h-negligible})}
and {\rm (\ref{eq:h-negligible-2})} we have
$$
\sqrt{\frac{nh}{\kappa_1}}\sqrt{\hat f(x)}\left(\widehat
{m^*}(x)-m_n^*(x)\right)\convdistr N(0,1) .
$$
\end{prop}
\begin{rem}{\rm
The result under (P1) follows by taking $\alpha_X=1$ in
(\ref{eq:h-negligible-2}). In this case the condition is always
fulfilled. The difference between the estimators of $m(\cdot)$ and
$m^*(\cdot)$ appears by comparison of (\ref{eq:h-negligible-1}) and
(\ref{eq:h-negligible}). In case of $m^*(\cdot)$ much larger
bandwidths are allowed to achieve classical rates of convergence.
The additional condition (\ref{eq:negligibility-condition}),
required for $m^*(\cdot)$ estimation is typically easy to verify,
even for long memory or non-stationary
sequences. \\

Recall from the previous remark that LRD in predictors basically does not matter. This is not the case for the
shape estimation. Assume that $hn^{1-\alpha_X}\to \infty$, but
{\rm (\ref{eq:negligibility-condition})}, {\rm (\ref{eq:h-negligible})}, (\ref{eq:h-negligible-2a}) and
\begin{equation}\label{eq:h-negligible-3}
\frac{\sqrt{nh}}{n^{\alpha_X/2}}\frac{1}{n}\sum_{i=1}^n\Exp[\varepsilon_i|{\cal H}_{i-1}]=o_P(1)
\end{equation}
hold. For example, this happens when $\varepsilon_i$, and thus
$\Exp[\varepsilon_i|{\cal H}_{i-1}]$ is the linear process with
$\Var\left(\sum_{i=1}^n\Exp[\varepsilon_i|{\cal H}_{i-1}]\right)\sim
C_1^2n^{2-\alpha}$, $\alpha\in (0,1)$, and $h^5n^{1-\alpha}\to 0$.
It follows from the proof of Proposition \ref{thm:main-1} that
$$
n^{\alpha_X/2}\left(\widehat
{m^*}(x)-m_n^*(x)\right)\convdistr {\cal N}(0,A_1^2\Exp^2[m(X_1)X_1]).
$$
}
\end{rem}

\subsection{Mean Square Error}\label{sec:mse}
In this section we establish asymptotic formulas for mean integrated squared error for both $m(\cdot)$ and $m^*(\cdot)$. In particular, it will be shown that we may improve the rates of convergence, if we estimate the shape function instead of $m(\cdot)$.  \\

Consider the following weighted versions of the mean integrated
squared errors:
\begin{eqnarray*}
\MISE_r(h):=\int\Exp\left[(\hat m_h(x)-m(x))^2\right]r(x)\;dx, \\
\MISE_r^*(h):=\int\Exp\left[(\hat m_h^*(x)-m^*(x))^2\right]r(x)\;dx,
\end{eqnarray*}
where $r(\cdot)$ is a weight function and integrals are taken over a support of $f$.
\begin{prop}\label{prop:mse_m}
Assume {\rm (P2)} and {\rm (E2)}. Then we have
\begin{eqnarray*}
\MISE_r(h)&\sim&\frac{\kappa_1}{nh}\int \frac{r(x)}{f(x)}\; dx+\frac{h^4\kappa_2^2}{4}\int  \left(\frac{m''(x)f(x)+2m'(x)f'(x)}{f(x)}\right)^2r(x)\; dx\\
&&+C_1^2n^{-\alpha}+C_1^2h^2\kappa_2n^{-\alpha}\int \frac{f''(x)}{f(x)}r(x)\; dx.
\end{eqnarray*}
\end{prop}
\begin{rem}{\rm
The first term in $\MISE_r(h)$ describes i.i.d. type behaviour, the second one is due to bias. The terms involving $n^{-\alpha}$ describe a possible contribution of long memory. Note that we have to include the term $h^2n^{-\alpha}$. These terms do not have influence on the optimal behaviour of $\MISE$, but they influence $h_{\rm opt}$, the optimal bandwidth choice. Indeed, one can verify that
$$
h_{\rm opt}\sim \left\{
\begin{array}{ll}
Cn^{-1/5} & \mbox{\rm if } \alpha>2/5;\\
Cn^{-(1-\alpha)/3} & \mbox{\rm if } \alpha<2/5,
\end{array}
\right.
$$
so that $\MISE_r(h_{\rm opt})$ is proportional to $n^{-4/5}$ if
$\alpha>4/5$ and $n^{-\alpha}$, if $\alpha<4/5$. Note also that
there is no contribution of LRD of the predictors. }
\end{rem}
\begin{prop}\label{lem:mse_hatm_ast}
Assume that
\begin{equation}\label{eq:var-assumption}
\Var\left(\sum_{i=1}^n(\Exp[\ve_{i} |
\mathcal{H}_{i-1}]-\ve_i)\right)=O(n)
\end{equation}
holds.
Under {\rm (P2)} and {\rm (E2)} we have
\begin{eqnarray*}
\lefteqn{\MISE_r^*(h)\sim\frac{\kappa_1}{nh}\int \frac{r(x)}{f(x)}\; dx}\\
&&+\frac{h^4}{2}\int  \left(\frac{m''(x)f(x)+2m'(x)f'(x)}{f(x)}\right)^2r(x)\; dx+A_1^2\Exp^2[m(X_1)X_1]n^{-\alpha_X}.
\end{eqnarray*}
\end{prop}
\begin{rem}{\rm
The condition (\ref{eq:var-assumption}) can be verified for many
time series, including time series with long memory (see Section
\ref{sec:examples}). The result under (P1) can be obtained by taking
$\alpha_X=1$; then the last term is negligible. Consequently, under
(P1) we remove long memory of errors, however, under (P2) there is
an additional term due to long memory of predictors.

Under (P1) the optimal bandwidth, $h_{\rm opt}^*$, is proportional
to $n^{-1/5}$, yielding $\MISE_r^*(h_{\rm opt}^*)\sim Cn^{-4/5}$. }
\end{rem}
\begin{rem}{\rm
Hall and Hart, \cite{HallHart1990}, were the first who proved the mean squared error behaviour in case of  fixed-design regression. The meaning of their results is that LRD in errors always influences estimation of the conditional mean.

On the other hand, in case of kernel density estimation based on LRD data $\varepsilon_1,\ldots,\varepsilon_n$, Hall and Hart \cite{HallHart1990a} showed a similar dichotomous behaviour, as described in Proposition \ref{prop:mse_m}.
}
\end{rem}
\subsection{Empirical bandwidth choice}\label{sec:bandwidth}
In this section we study properties of empirical bandwidth selector procedures.
We shall focus on the case of i.i.d. predictors, to show an influence of LRD errors on the empirical risk. We will work also under the two additional assumptions: $\Exp[m^2(X_1)]<\infty$ and that $f$ has bounded support. This will simplify some computations and allow us to write $\MISE_f(h)$.
Let
$$
\ASE(h)=\frac{1}{n}\sum_{i=1}^n\left(\hat m_{h}(X_i)-m(X_i)\right)^2
$$
be the Averaged Squared Error.\\

First, we answer the question, whether minimization of $\ASE$ leads to a valid minimizer. The answer is affirmative: the meaning of (\ref{eq:conclusion-1}) is that the quotient of $\hat h$, the minimizer of $\ASE$, and the minimizer of $ \MISE_{f}$ converges to 1 in probability. Furthermore, $\ASE(\hat h)/\MISE(h_{\rm opt})\convprob 1$.
\begin{prop}\label{prop-ISE}
Assume that $f, m \in {\cal C}^2$.
Let $B_1<B_2$ be finite and positive constants. Under {\rm (P1)} and {\rm (E2)} we have
uniformly over $[B_1h_{\rm opt},B_2h_{\rm opt}]$,
\begin{eqnarray}\label{eq:conclusion-1}
\ASE(h)- \MISE_{f}(h)&=&o_P\left( \MISE_{f}(h_{\rm opt})\right).
\end{eqnarray}
\end{prop}
However, what we are interested in from a practical point of view, is, first, whether a cross-validation produces a valid bandwidth, and, second,  whether a cross-validation is a good approximation to $\ASE$ and $\MISE$.
To answer this, let $\hat m_{j,h}(\cdot)$ be the version of the estimator (\ref{eq:kernel-estim}), where the summation is over $i\not\in I_j(l)$, where $I_j:=\{i:|j-i|>l\}$. Empirical cross-validation bandwidth is obtained via minimizing
$$
\CV_l(h) :=\frac{1}{n}\sum_{i=1}^n (Y_i-\hat m_{i,h}(X_i))^2.
$$
Denote $\CV(h)=\CV_0(h)$. Note that if both predictors and errors are i.i.d., then
$$
\Exp [\CV(h)]=\Exp[\ASE(h)]+\Exp[\varepsilon_1^2]=\MISE_{f}(h)+\Exp[\varepsilon_1^2],
$$
i.e. in the average sense $\CV(h)$ is the exact approximation of $\MISE_{f}(h)+\Exp[\varepsilon_1^2]$.\\

The result for LRD data is as follows.
\begin{prop}\label{prop-CV}
Assume that $f, m \in {\cal C}^2$.
Let $B_1<B_2$ be finite and positive constants. Under {\rm (P1)} and {\rm (E2)} we have uniformly over
$[B_1h_{\rm opt},B_2h_{\rm opt}]$,
\begin{equation*}
\frac{\CV(h)- \MISE_{f}(h)-\frac{1}{n}\sum_{i=1}^n\varepsilon_i^2}{ \MISE_{f}(h_{{\rm opt}})}=\frac{1}{ \MISE_{f}(h_{{\rm opt}})}\frac{1}{n^2}\sum_{j,j'=1\atop j\not=j'}^n\varepsilon_j\varepsilon_{j'}+o_p(1).
\end{equation*}
\end{prop}
Let us comment on the above result.
For a fixed-design case, under appropriate conditions on $l$, we have (see \cite{HallLahiriPolzehl1995})
\begin{equation}\label{eq:Hall-1}
\CV_l'(h)-\MISE'(h)-\frac{1}{n}\sum_{i=1}^n\varepsilon_i^2=o_P(n^{-4\alpha/5}),
\end{equation}
uniformly over $[Cn^{-\alpha/5},C'n^{-\alpha/5}]$, $C<C'$. (Here,
$\CV_l'(h)$ and $\MISE'(h)$ are defined in a slightly different way,
to accommodate fixed-design). Note that $n^{-4\alpha/5}$ is the rate
of $\MISE'(h_{{\rm opt}}')$ and $n^{-\alpha/5}$ is asymptotically
proportional to $h_{{\rm opt}}'$, where the latter is the
asymptotically optimal bandwidth in the fixed design regression.
This means that the ratio of the bandwidth obtained by
cross-validation and the $\MISE$ optimal bandwidth converges to 1 in
probability. Also, $\Exp[\CV_l'(h_{\rm opt}')]$ provides a valid
approximation to $\MISE'(h_{{\rm opt}}')+\Exp[\varepsilon_1^2]$.
(The last statement is intuitive only, since the rate at the right
hand side of (\ref{eq:Hall-1}) is in probability rather than in
$L^1$).

Now, from Proposition \ref{prop-CV} we conclude that in the
random-design regression $\hat h_{\rm CV}$, the minimizer of
$\CV(h)$, has the property $\hat h_{\rm CV}/h_{\rm opt}\convprob 1$.
However, $\CV(h)$ itself provides a valid approximation only if
$\alpha>4/5$. This agrees with the results in
\cite{HallLahiriTruong1995} in case of density estimation.
\section{Numerical studies}\label{sec:numerical}
Simulation studies were conducted as follows:
\begin{enumerate}
\item We set $h$.
\item Simulate 100 errors $\varepsilon_i$ from \texttt{ARIMA} with different LRD parameters $d$. Note that $d=(1-\alpha)/2$ and i.i.d.~case corresponds to $d=0$. We used \texttt{R}-package \texttt{fracdiff}.
\item Simulate 100 predictors following standard normal random variables. First, we simulate i.i.d. predictors, then LRD predictors with $d_X:=(1-\alpha_X)/2=0.3$.
\item This procedure was repeated 500 times.
\item As the output we get a  Monte Carlo approximation to the $\MISE$ and $\MISE^*$.
\end{enumerate}
Table \ref{tab:1} shows results for both $m(\cdot)$ and $m^*(\cdot)$
for the function $m(x)=\sin(2\pi x)$ and bandwidths $h=0.05$, $h=1$,
respectively. Note that in this case $m=m^*$. Even for a relatively
small sample size, we may observe that $\MISE^*$ for the shape
function remains constant for either choice of the bandwidth
($h=0.05$ or $h=1$). On the other hand, for the small bandwidth
$h=0.05$ we observe that $\MISE$ for $m(\cdot)$ stays constant up to
$d=0.25$, but it grows almost immediately for $h=1$. Furthermore,
LRD starts to dominate earlier for the larger bandwidth. This is in
line with Propositions \ref{prop:mse_m} and \ref{lem:mse_hatm_ast}.
\begin{table}[h]
\begin{center}
\begin{tabular}{||c||c|c|c|c||}\hline\hline
$d$ & \multicolumn{2}{c|}{$h=0.05$} &  \multicolumn{2}{c||}{$h=1$}     \\\hline
    & $\MISE$ & $\MISE^*$ & $\MISE$ & $\MISE^*$\\\hline
0       & 0.3995014 & 0.3983940 & 0.5087492 &0.5016983\\
0.05     & 0.3932386 & 0.3870820 & 0.5104933 &0.5011370 \\
0.10    & 0.3859548 & 0.3791453 & 0.5136161 &0.5042601\\
0.15     & 0.3780757 & 0.3724782 & 0.5163214 &0.5030991\\
0.20    & 0.4095343 & 0.3910381 & 0.5255956&0.5026314 \\
0.25     & 0.4050078 & 0.3850845 &0.5322871 &0.5027839\\
0.30    & 0.4242228 & 0.3794874 & 0.5494079& 0.5019103\\
0.35    & 0.4493550 & 0.3940342 & 0.5975138&0.5044394\\
0.40     & 0.5848862 & 0.3805059 & 0.6432639&0.5016800\\
0.45    & 0.8899494 & 0.3822220 & 0.8908165&0.5038973\\\hline\hline
\end{tabular}
\end{center}
\caption{ {\small  MISE for some values
of the dependence parameter $d$ and i.i.d. predictors.} }\label{tab:1}
\end{table}
Next (Table \ref{tab:2}), we repeated this experiment with dependent predictors, choosing $d_X:=(1-\alpha_X)/2=0.3$. By comparing both tables, note that there is a little influence of LRD of predictors on $\MISE$ for $m(\cdot)$. This is still in line with the Proposition \ref{prop:mse_m}. On the other hand, for $m^*(\cdot)$ estimation, Proposition \ref{lem:mse_hatm_ast} suggests that LRD of predictors should contribute. It does not seem to be the case here, however, or simulation studies suggest that $\MISE^*$ depends on the constant $\Exp^2[m(X_1)X_1]$, as indicated by our theoretical results.

\begin{table}[h]
\begin{center}
\begin{tabular}{||c||cc||}\hline\hline
$d$& \multicolumn{2}{c||}{$h=0.05$}   \\\hline
    & $\MISE$ & $\MISE^*$ \\\hline
0     & 0.4021701&0.3963913  \\
0.05  &0.4454991&0.4385924  \\
0.10    &0.4337356 &0.4227792 \\
0.15    &0.4281452  &0.4225283\\
0.20    &0.4262649  &0.4121518 \\
0.25    &0.4436659  &0.4137239 \\
0.30   & 0.4417674 &0.4049835\\
0.35   &0.4948932  &0.4307672\\
0.40   &0.5673097&0.4146042  \\
0.45  & 0.8217908&0.4258788\\\hline\hline
\end{tabular}
\caption{ {\small MISE for some values
of the dependence parameter $d$ and LRD predictors.}}\label{tab:2}
\end{center}
\end{table}

Finally, based on 500 simulations, we computed averaged values of $\CV$ criterion with optimally chosen $h$. Table \ref{tab:3} indicates that LRD influences $\CV$ almost immediately, which is once again in line with our theoretical results.
\begin{table}[h]
\begin{center}
\begin{tabular}{||c||c||}\hline\hline
d    & $\CV$  \\\hline
0     & 0.2912500 \\
0.05  & 0.2921772\\
0.10    & 0.2960790 \\
0.15    & 0.2972158 \\
0.20    & 0.310086\\
0.25    & 0.3206689 \\
0.30   &  0.3415160\\
0.35   & 0.3687273\\
0.40   &  0.4133771\\
0.45  & 0.4663486\\\hline\hline
\end{tabular}
\caption{ {\small $\CV$ for some values
of the dependence parameter $d$ and i.i.d. predictors.}}\label{tab:3}
\end{center}
\end{table}

\section{Examples}\label{sec:examples}
In this section we present some examples, which show that our conditions are easy to verify for very different
long memory processes. \\

Unless specified otherwise, $\{Z,Z_i,i\ge
1\}$ and $\{\eta,\eta_i,-\infty<i<\infty\}$ will be sequences of
centered i.i.d. random variables with
$\Exp[Z_1^1]=\Exp[\eta_1^2]=1$. More detailed description of most of
the models below, together with stationarity and moment conditions,
can be found in \cite{GirLeiSur2007}.
\begin{ex}[Linear processes]\label{ex:linear}{\rm
Assume (E2). Then we have
\begin{equation}\label{eq:1}
\Exp[\varepsilon_i|{\cal H}_{i-1}]=\sum_{k=1}^{\infty}c_k\eta_{i-k} +\Exp[\eta_i|{\cal H}_{i-1}]=:\varepsilon_{i,i-1}+0
\end{equation}
and $\{\varepsilon_{i,i-1},i\ge 1\}$ is LRD linear process with (up to a constant) the same limiting behavior as $\{\varepsilon_i,i\ge 1\}$.
Consequently,
\begin{equation}\label{ex:1}
\sum_{i=1}^n\left\{\Exp[\varepsilon_i|{\cal
H}_{i-1}]-\varepsilon_i\right\}=-\sum_{i=1}^n\eta_i,
\end{equation}
and
\begin{equation}\label{ex:2}
\sum_{i=1}^n\Exp[\varepsilon_i|{\cal H}_{i-1}]=\sum_{i=1}^n\varepsilon_{i,i-1}.
\end{equation}
By (\ref{ex:1}) and since $\eta_i$ are i.i.d.,
(\ref{eq:negligibility-condition}) is automatically fulfilled.
Finally, by (\ref{ex:2}), we conclude that (\ref{eq:h-negligible})
holds if
\begin{equation}\label{eq:small-h}
h^5n^{1-\alpha}\to 0.
\end{equation}
Note that this condition is much less restrictive than
\begin{equation}\label{eq:small-h-usual}
hn^{1-\alpha}\to 0,
\end{equation}
which is required for (\ref{eq:h-negligible-1}) to hold. The latter condition is the same as in \cite{WuMielniczuk2002}.
}\end{ex}
\begin{ex}[Functionals of Linear Processes]\label{ex:functional-of-linear}{\rm
Consider the linear process from Example \ref{ex:linear}. Let $T$ be
a twice differentiable functional and let
\begin{equation*}
\varepsilon_i=T\left(\sum_{k=0}^{\infty}c_k\eta_{i-k}\right)=T(c_0\eta_i+\varepsilon_{i,i-1}).
\end{equation*}
Consider the same assumptions as in Example \ref{ex:linear}.
Additionally, we assume that $\Exp[T(\varepsilon_1)]=0$.

Let $f_{\eta}$ be the density of $\eta_1$. Then, by considering two terms of Taylor expansion (which is enough for a quadratic functional),
\begin{eqnarray*}
\lefteqn{\Exp[\varepsilon_i|{\cal H}_{i-1}]= \int T(u+\varepsilon_{i,i-1})f_{\eta}(u)du}\\
&=& \int T(u)f_{\eta}(u)du+\varepsilon_{i,i-1} \int T'(c_0u)f_{\eta}(u)du+\frac{1}{2}\varepsilon_{i,i-1}^2\int T^{''}(u)f_{\eta}(u)du\\
&=& \Exp[T(\eta_1)]+\Exp[T'(\eta_1)]\varepsilon_{i,i-1}+\frac{1}{2}\Exp[T^{''}(\eta_1)]\varepsilon_{i,i-1}^2
\end{eqnarray*}
and
$$
\varepsilon_i=T(\eta_i)+T'(\eta_i)\varepsilon_{i,i-1}+\frac{1}{2}T^{''}(\eta_i)\varepsilon_{i,i-1}^2.
$$
For simplicity take $T(u)=u^2-\Exp[\eta_1^2+\varepsilon_{1,0}^2]$ and assume
that the density $f_{\eta}$ is symmetric. Then
$\Exp[T(\eta_1)]=\Exp[\varepsilon_{i,i-1}^2]=-\sum_{k=1}^{\infty}c_k^2$,
$\Exp[T'(c\eta_1)]=2\Exp[\eta_1]=0$. Consequently,
$$
\Var\left(\sum_{i=1}^n\Exp[\varepsilon_i|{\cal H}_{i-1}]\right)=\Var\left(\sum_{i=1}^n\left(\varepsilon_{i,i-1}^2 -\Exp[\varepsilon_{i,i-1}^2]\right)\right)\sim
O(n^{2-2\alpha}\vee n).
$$
We conclude that (\ref{eq:h-negligible}) holds for all $\alpha>1/2$ or if $n^{1-2\alpha}h^5\to 0$, whereas (\ref{eq:h-negligible-1}) holds when $\alpha>1/2$ or if $n^{1-2\alpha}h\to 0$.

Moreover,
$$
\Exp[\varepsilon_i|{\cal H}_{i-1}]-\varepsilon_i=-\eta_i-2\eta_i\varepsilon_{i,i-1}
$$
and thus (\ref{eq:negligibility-condition}) holds since $\eta_i,i\ge 1$ are i.i.d. and $\eta_i\varepsilon_{i,i-1},i\ge 1$ are uncorrelated.

Similar consideration can be carried out for any functional $T$, in
particular, for $T(u)=|u|^{\delta}-\Exp[|\varepsilon_1|^{\delta}]$. The
set of parameters for which (\ref{eq:h-negligible}) and
(\ref{eq:h-negligible-1}) hold depends on the so-called power rank
of $T$ (see \cite{HoHsing1997} for more details). If the power rank
is 1, the (\ref{eq:h-negligible}) and (\ref{eq:h-negligible-1}) hold
for $\alpha$, $h$ as in Example \ref{ex:linear}, if the power rank is 2,
then the conditions are fulfilled for $\alpha$, $h$ as in case of
quadratic functional discussed above. }
\end{ex}
\begin{ex}[FARIMA-GARCH processes]\label{ex:farima-garch}{\rm
Assume that
$$
\varepsilon_i=(1-B)^{-d}\phi^{-1}(B)\psi(B)\eta_i,
$$
where $\eta_i=Z_ih_i^{1/2}$ and $h_i$ is GARCH($r,s$),
$$
h_i=a_0+\sum_{j=1}^ra_j\eta_{i-j}^2+\sum_{k=1}^s\beta_kh_{i-k}.
$$
Here, $B$ is the backshift operator, $\psi$ and $\phi$ are polynomials in $B$ and $d\in (-1/2,1/2)$. Note that under appropriate stationarity conditions the FARIMA-GARCH process can be written as the linear process in (E2), where $c_k\sim Ck^{-\beta}$, $\beta=1-d$ (we refer to \cite{BeranFeng2001} for more details).

Let ${\cal H}_{i-1}=\sigma(\eta_i,Z_i,\eta_{i-1},Z_{i-1}\ldots)$.
Then $\Exp[\eta_i|{\cal H}_{i-1}]=\Exp[Z_i]\Exp[h_i^{1/2}|{\cal
H}_{i-1}]=0$. Consequently, as in Example \ref{ex:linear} (see
(\ref{eq:1})), (\ref{eq:negligibility-condition}) holds for all
$d\in (-1/2,1/2)$, since $\eta_i$ are uncorrelated. Furthermore,
(\ref{eq:h-negligible}) and (\ref{eq:h-negligible-1}) hold if
$n^{2d}h^5\to 0$ and $n^{2d}h\to 0$, respectively, on account of
\cite[Theorem 3]{BeranFeng2001}. }
\end{ex}
\begin{ex}[Antipersistent errors]{\rm
If in Example \ref{ex:farima-garch} $d\in (-1/2,0)$, then we have
{\it antipersistence} and (\ref{eq:h-negligible}) is always
fulfilled. Consequently, in case of antipersistent errors the
correct scaling for the estimator of both $m$ and $m^*$ is always
$\sqrt{nh}$. Note that in case of fixed-design regression
antipersistency may improve convergence beyond i.i.d. rates, see
e.g. \cite{BeranFeng2001}, \cite{BeranFeng2002}. }
\end{ex}
\begin{ex}[Stochastic volatility]\label{ex:sv}{\rm
Let $T$ and $W$ be real-valued functionals.
Define $$\varepsilon_i^*=Z_iR_i, \qquad R_i=W\left(a+\sum_{k=1}^{\infty}c_k\eta_{i-k}\right),\quad a>0.$$
Let $$\varepsilon_i=T(\varepsilon_i^*)-\Exp[T(\varepsilon_i^*)]$$
and ${\cal H}_i=\sigma(\eta_i,Z_i,\eta_{i-1},Z_{i-1},\ldots)$. Note that $R_i$ is ${\cal H}_{i-1}$-measurable.
For simplicity, assume that $T(uv)=T(u)T(v)$, which applies e.g. to polynomials $T(u)=|u|^{\delta}$.
Then
\begin{eqnarray*}
\sum_{i=1}^n \Exp[\varepsilon_i|{\cal H}_{i-1}]&=&\Exp[T(Z)]\sum_{i=1}^n\left\{\Exp[T(R_i)|{\cal H}_{i-1}]-\Exp[T(R_i)]\right\}\\
&=& \Exp[T(Z)]\sum_{i=1}^n\left\{T(R_i)-\Exp[T(R_i)]\right\}.
\end{eqnarray*}
Thus, if $c_k\sim k^{-(\alpha+1)/2}$, $\alpha\in (0,1)$, the
conditions for (\ref{eq:h-negligible}) and (\ref{eq:h-negligible-1})
are the same as for nonlinear transformations of linear processes in
Example \ref{ex:functional-of-linear} by substituting $T\to T\circ
W$.

Furthermore
$$
\sum_{i=1}^n\left\{\Exp[\varepsilon_i|{\cal H}_{i-1}]-\varepsilon_i\right\}=\sum_{i=1}^nT(R_i)\left\{\Exp[T(Z)]-T(Z_i)\right\}=:\sum_{i=1}^nT(R_i)U_i.
$$
Note that $T(R_i)U_i,i\ge 1$, are uncorrelated, thus (\ref{eq:negligibility-condition}) is always fulfilled. \\

If $\varepsilon_i=\varepsilon_i^*=Z_iR_i$, then $\Exp[\varepsilon_i|{\cal H}_{i-1}]=0$ and since the random variables $\varepsilon_i,i\ge 1$ are uncorrelated, $\sum_{i=1}^n\varepsilon_i=O_P(\sqrt{n})$ . Thus, the memory parameter $\alpha$ has no influence on the asymptotic behavior of neither $\hat m_h$ nor $\hat m_{h}^*$.
}
\end{ex}
\begin{ex}[LARCH processes]{\rm
Define $$\varepsilon_i^*=Z_iR_i, \qquad R_i=a+\sum_{k=1}^{\infty}c_k\varepsilon_{i-k}^*,\quad a>0$$
and assume that $\sum_{k=1}^{\infty}c_k^2<1$, $c_k\sim k^{-(\alpha+1)/2}$, $\alpha\in (0,1)$.
Let $$\varepsilon_i=T(\varepsilon_i^*)-\Exp[T(\varepsilon_i^*)]$$
and ${\cal H}_i=\sigma(\varepsilon_i^*,Z_i,\epsilon_{i-1}^*,Z_{i-1},\ldots)$. The random variable $R_i$ is ${\cal H}_{i-1}$-measurable. As in Example \ref{ex:sv}, assume that $T(uv)=T(u)T(v)$ so that
$$
\sum_{i=1}^n\left\{\Exp[\varepsilon_i|{\cal H}_{i-1}]-\varepsilon_i\right\}=\sum_{i=1}^nT(R_i)\left\{\Exp[T(Z)]-T(Z_i)\right\}=:\sum_{i=1}^nT(R_i)U_i.
$$
and $T(R_i)U_i,i\ge 1$ are uncorrelated, thus (\ref{eq:negligibility-condition}) is always fulfilled. Furthermore,
$$
\sum_{i=1}^n \Exp[\varepsilon_i|{\cal H}_{i-1}]=\Exp[T(Z)]\sum_{i=1}^n\left\{T(R_i)-\Exp[T(R_i)]\right\}.
$$
Although this expression has the same form as in Example
\ref{ex:sv}, $R_i$ is not a linear process based on i.i.d. random
variables. Nevertheless, from \cite[Theorem 1.1]{BerHor2003} we
conclude that if $T$ is twice differentiable, and
$\Exp[R_1T(R_1)]\not=0$, then the scaling factor for the latter sum
is the same as for linear processes in Example \ref{ex:linear}.
Thus, for the conditions (\ref{eq:h-negligible}) and
(\ref{eq:h-negligible-1}) to be fulfilled, we need
(\ref{eq:small-h}) and (\ref{eq:small-h-usual}), respectively. }
\end{ex}
\section{Proofs}\label{sec:proofs}
In the proofs we apply a concept of martingale approximation and Hermite expansion. In the
context of nonparametric estimation the first method was introduced in
\cite{WuMielniczuk2002} and \cite{MielniczukWu2004}, the latter one is a standard tool in LRD setting, see e.g. \cite{Taqqu2003}.

Let us note that under the regularity assumptions we have:
\begin{eqnarray}
 \Exp[K_h(x-X_1)]& = & hf(x)+h^3 \cdot {f''(x)\over 2} \kappa_2+o(h^3),  \label{eq:EKh}
\end{eqnarray}
We may write
\begin{equation}\label{eq:decomposition}
\hat m_h(x)-m_n(x)= \frac{\hat f_h(x)}{f(x)}\left(\hat m_h(x)-m_n(x)\right)+(\hat m_h(x)-m_n(x))\frac{(f(x)-\hat f_h(x))}{f(x)}.
\end{equation}
Since $\hat f_h$ is the consistent estimator of $f$, it suffices to study the first part only.
Decompose
\begin{eqnarray*}
\lefteqn{\frac{\hat f_h(x)}{f(x)}\left(\hat m_h(x)-m_n(x)\right) = \frac{1}{nhf(x)}\sum_{i=1}^nK_h\left(x-X_i\right)\left(m(X_i)-m_n(x)\right)}\nonumber\\
&&+
\frac{1}{nhf(x)}\sum_{i=1}^nK_h\left(x-X_i\right)\ve_i=:N_n(x)
+ N_n'(x).\qquad\qquad\qquad\qquad
\end{eqnarray*}
The parts $N_n(\cdot)$ and $N_n'(\cdot)$ are decomposed further as follows:
\begin{eqnarray*}
\lefteqn{N_n(x)=\frac{1}{f(x)}(m(x)-m_n(x))(\hat f_h(x)-\Exp[\hat f_h(x)])}\\
&&+\frac{1}{nhf(x)}\sum_{i=1}^n\left(K_h\left(x-X_i\right)(m(X_i)-m(x))-\Exp[K_h\left(x-X_i\right)(m(X_i)-m(x))|{\cal X}_{i-1}]\right)\\
&&+
\frac{1}{nhf(x)}\sum_{i=1}^n\left(\Exp[K_h\left(x-X_i\right)(m(X_i)-m(x))|{\cal X}_{i-1}]-\Exp[K_h\left(x-X_i\right)(m(X_i)-m(x))]\right)\\
&=:&N_{n,1}(x)+N_{n,0}(x)+N_{n,2}(x).
\end{eqnarray*}
Likewise,
\begin{eqnarray}
\lefteqn{N_n'(x)=\frac{1}{nhf(x)}\sum_{i=1}^nK_h(x-X_i)\ve_i= }\nonumber\\&=&
\frac{1}{nhf(x)}\sum_{i=1}^n\left(K_h(x-X_i)\ve_i-\Exp\left[K_h(x-X_i)\ve_i|{\cal X}_{i-1}\vee {\cal H}_{i-1}\right]\right)\nonumber\\
&&+\frac{1}{nhf(x)}\sum_{i=1}^n\Exp\left[K_h(x-X_i)\ve_i|{\cal X}_{i-1}\vee {\cal H}_{i-1}\right]=:M_n(x)+D_n(x)\nonumber.
\end{eqnarray}
Consequently, for $\hat m_h(x)-m_n(x)$ we have the following
decomposition:
\begin{equation}\label{eq:function-LRD}
\frac{\hat f_h(x)}{f(x)}\left(\hat m_h(x)-m_n(x)\right)=N_{n,0}(x)+N_{n,1}(x)+N_{n,2}(x)+M_n(x)+D_n(x).
\end{equation}
Note that $M_n(\cdot)$ and $N_{n,0}(\cdot)$ are always martingales. \\

Assume first (P1). Then $N_{n,2}(\cdot)\equiv 0$. Furthermore,
\begin{equation}\label{eq:Dn-p1}
D_n(x)=\frac{1}{nhf(x)}\Exp[K_h(x-X_1)]\sum_{i=1}^n\Exp[\varepsilon_i|{\cal H}_{i-1}].
\end{equation}

Assume now (P2). Let $X_{i,i-1}=\Exp[X_i|{\cal X}_{i-1}]=\sum_{k=1}^{\infty}a_k\zeta_{i-k}$.
Let $f_{\zeta_i}(\cdot)$ be the density of $\zeta_i$. Let
$\gamma^2=\Var(X_{1,0})$, write $Z_i=X_{i,i-1}/\gamma$ and
note that $Z_i$ is independent of $\zeta_i$.
Let $H_q(\cdot)$ be the $q$th Hermite polynomial. Applying the
Hermite expansion we represent $N_{n,2}(x)$ as
\begin{eqnarray*}
N_{n,2}(x)=\frac{1}{nhf(x)}
\sum_{q=1}^{\infty}\frac{1}{q!}\sum_{i=1}^nH_q(Z_i) \int L(q;u)f_{\zeta_1}(u)\;
du,
\end{eqnarray*}
where
$$
L(q;u)=L(q;u,h,\gamma)=\Exp[K_h(x-(u+\gamma
Z_1))\left(m(u+\gamma Z_1)-m(x)\right)H_q(Z_1)].
$$
Also,
\begin{equation*}
D_n(x)=\frac{1}{nhf(x)}\sum_{q=0}^{\infty}\frac{1}{q!}\sum_{i=1}^nH_q(Z_i)\Exp[\varepsilon_i|{\cal H}_{i-1}]\int J(q;u)f_{\zeta_1}(u)\;
du,
\end{equation*}
where
$$
J(q;u)=J(q;u,h,\gamma)=\Exp[K_h(x-(u+\gamma
Z_1))H_q(Z_1)].
$$
Note that summations in Hermite expansions is from $q=0$, since the expanded functions does not have mean 0 w.r.t. Gaussian density.  \\

\noindent Furthermore,
\begin{equation}\label{eq:Nn1}
N_{n,1}(x)=\frac{1}{f(x)}(m(x)-m_n(x))\frac{1}{nh}\sum_{q=1}^{\infty}\sum_{i=1}^n\frac{C(q)}{q!}H_q(X_i),
\end{equation}
where
$$
C(q)=C(q;h)=\Exp[K_h(x-X_1)H_q(X_1)].
$$

As for the shape function $m^*$, we write
$$
\hat m_h^*(x)-m^*(x)=\left(\hat m_h^*(x)-m_n^*(x)\right) + \left(m_n^*(x)-m^*(x)\right),
$$
where $m_n^*(x)=m_n(x)-\int mf$. Clearly, $m_n^*(x)-m^*(x)=m_n(x)-m(x)$.
Let $\hat\Theta_1=\frac{1}{n}\sum_{i=1}^nm(X_i)$, $\hat\Theta_2=\frac{1}{n}\sum_{i=1}^n\varepsilon_i$. With this notation we write
\begin{eqnarray}
\lefteqn{\frac{\hat f_h(x)}{f(x)}\left(\hat m_h^*(x)-m_n^*(x)\right)=N_n(x)+\frac{1}{nhf(x)}\sum_{i=1}^nK_h(x-X_i)\left(\int mf-\hat\Theta_1\right)}\nonumber\\
&&+\frac{1}{nhf(x)}\sum_{i=1}^nK_h\left(x-X_i\right)\left(\ve_i-\hat\Theta_2\right)\nonumber\\
&=& N_{n,0}(x)+M_n(x)+\frac{1}{nhf(x)}\sum_{i=1}^nK_h(x-X_i)\left(\int mf-\hat\Theta_1\right)\nonumber\\
&&+D_n(x)-\frac{1}{nhf(x)}\sum_{i=1}^nK_h\left(x-X_i\right)\hat\Theta_2+N_{n,1}(x)+N_{n,2}(x).\label{eq:star-expansion}
\end{eqnarray}
The crucial difference between $\hat m_h(\cdot)$ and $\hat m_h^*(\cdot)$ is that possibly LRD part $D_n(x)$ is replaced with
$$D_n(x)-\frac{1}{nhf(x)}\sum_{i=1}^nK_h\left(x-X_i\right)\hat\Theta_2=:D_n(x)-E_n(x).$$
\subsection{Proof of Propositions \ref{thm:main}, \ref{thm:main-1}}
Recall (\ref{eq:function-LRD}).
Since $N_{n,0}(x)$ and $M_n(x)$ are martingales we may easily conclude that
\begin{equation}\label{eq:9}
\Var\left( N_{n,0}(x)\right)=O(h/n),\quad \Var \left(M_n(x)\right)\sim \frac{\kappa_1}{nhf(x)},
\end{equation}
which means that $N_{n,0}(x)$ is negligible. Furthermore, the martingale part $M_n(x)$ may be studied using standard tools (see the proof below in Section \ref{sec:proof-mtg}).
\begin{lem}\label{lem:mtg-convergence}
Under the conditions of Proposition \ref{thm:main},
\begin{equation*}
\sqrt{\frac{nh}{\kappa_1}}\sqrt{f(x)}M_n(x)\convdistr N(0,1).
\end{equation*}
\end{lem}
Under (P1), using (\ref{eq:bias}),
$$
\Var\left[N_{n,1}(x)\right]=O(h^4(nh)^{-1})=O(h^3/n),
$$
so that this term is negligible w.r.t. $M_n(\cdot)$ as well.
Furthermore, using (\ref{eq:Dn-p1}) and (\ref{eq:EKh}) we have
\begin{equation}\label{eq:10}
D_n(x)=\left(1+h^2\frac{\kappa_2}{2}\frac{f^{''}(x)}{f(x)}+o(h^2)\right)\frac{1}{n}\sum_{i=1}^n\Exp[\varepsilon_i|{\cal
H}_{i-1}].
\end{equation}
Consequently, under (P1) the result follows form Lemma \ref{lem:mtg-convergence} and assumption (\ref{eq:h-negligible-1}), which makes $D_n(\cdot)$ negligible.\\

Now, we work under the assumption (P2). Recall (\ref{eq:Nn1}). We split the stochastic term there as
$$
C(1)\sum_{i=1}^nX_i+\sum_{q=2}^{\infty}\sum_{i=1}^n\frac{C(q)}{q!}H_q(X_i).
$$
Using orthonormality of the Hermite polynomials, the variance of the latter term is
$$
\sum_{i=1}^n\sum_{q=2}^{\infty}\frac{C^2(q)}{q!}\Cov^q(X_1,X_i)\le \sum_{i=1}^n\Cov^2(X_1,X_i)\|K_h(x-\cdot)\|,
$$
where $\|\cdot\|$ stands for $L^2$ norm with respect to the Gaussian measure. Since $\|K_h(x-\cdot)\|$ and
$C(1)\sim h xf(x)$, we conclude that the leading term in $N_{n,1}(x)$ is
$$
\frac{1}{f(x)}(m(x)-m_n(x))\frac{C(1)}{hf(x)}\frac{1}{n}\sum_{i=1}^nX_i.
$$
This implies
\begin{equation}\label{eq:7a}
\Var (N_{n,1}(x))=O(h^4n^{-\alpha_X}).
\end{equation}
The similar consideration is applied to $D_n$ and $N_{n,2}$. Using the independence of
$Z_1$ and $\zeta_1$, we have
\begin{equation}\label{eq:8}
\int J(0;u)f_{\zeta_1}(u)\;
du=\Exp[K_h(x-X_1)]=hf(x)+\frac{1}{2}h^3f''(x)\int u^2K(u)\; du+o(h^3).
\end{equation}
This yields
\begin{equation}\label{eq:7}
D_n(x)=\frac{1}{n}\sum_{i=1}^n\Exp[\varepsilon_i|{\cal H}_{i-1}]+o_P(1).
\end{equation}
The leading terms in $N_{n,2}(x)$ is
$$
\frac{\Exp[L(1;\zeta_1)]}{hf(x)}\frac{1}{n}\sum_{i=1}^nZ_i,
$$
which implies
\begin{equation}\label{eq:7b}
\Var (N_{n,2}(x))=O(h^4n^{-\alpha_X}).
\end{equation}
Consequently, under (P2), the result follows from Lemma \ref{lem:mtg-convergence} and assumption (\ref{eq:h-negligible-1}), which makes $D_n(\cdot)$ negligible, together with (\ref{eq:h-negligible-2a}), which makes $N_{n,1}(x)+N_{n,2}(\cdot)$ negligible. \\

We prove now Proposition \ref{thm:main-1} assuming (P2). Recall (\ref{eq:star-expansion}). It was proven before that $N_{n,0}(x)$, $N_{n,1}(x)$ and $N_{n,2}(x)$ are negligible. The first term in the Hermite expansion for $N_{n,2}'(x)$ is
$$\frac{1}{nhf(x)}\Exp[J(0;\zeta_1)]\sum_{i=1}^n\Exp[\varepsilon_i|{\cal H}_{i-1}].
$$
If $\alpha+\alpha_X>1$, then $\sum_{i=1}^n\Exp[\varepsilon_i|{\cal H}_{i-1}]Z_i=O_P(\sqrt{n})$, otherwise
$\sum_{i=1}^n\Exp[\varepsilon_i|{\cal H}_{i-1}]Z_i=O_P(n^{1-(\alpha+\alpha_X)/2})$.
Since for $q\ge 0$, $\Exp[J(q;\zeta_1)]=O(h)$, we conclude that $D_n(x)$ can be written as
$$
\frac{1}{nhf(x)}\Exp[J(0;\zeta_1)]\sum_{i=1}^n\Exp[\varepsilon_i|{\cal H}_{i-1}]+O_P(n^{-(\alpha+\alpha_X)/2})+O_P(n^{-1/2}).
$$
The first two terms in the Hermite expansion for $E_n(x)$ are
$$
\Theta_2\frac{1}{hf(x)}\Exp[K_h(x-X_1)]+\Theta_2\frac{1}{hf(x)}\Exp[K_h(x-X_1)X_1]\frac{1}{n}\sum_{i=1}^nX_i.
$$
Using (\ref{eq:8}), and $\frac{1}{n}\sum_{i=1}^nX_i=O_P(n^{-\alpha_X/2})$,
we conclude that the leading terms in the difference $N_{n,2}'(x)-E_n(x)$, are
\begin{eqnarray}\label{eq:Nn2}
\lefteqn{\frac{1}{n}\sum_{i=1}^n\left(\Exp[\varepsilon_i|{\cal H}_{i-1}]-\varepsilon_i\right)}\nonumber\\
&&+O_P\left(\frac{h^2}{n}\sum_{i=1}^n\Exp[\varepsilon_i|{\cal H}_{i-1}]\right)+
O_P\left(n^{-\alpha_X/2}\frac{1}{n}\sum_{i=1}^n\Exp[\varepsilon_i|{\cal H}_{i-1}]\right).
\end{eqnarray}
The first two terms are negligible under the conditions (\ref{eq:negligibility-condition}) and (\ref{eq:h-negligible}), respectively. The last term is negligible on account of (\ref{eq:h-negligible-3}), which is weaker than
(\ref{eq:h-negligible-2}). Finally,
$$
n^{\alpha_X/2}\left(\int mf-\hat\Theta_1\right)\convdistr {\cal N}(0,A_1^2\Exp^2[m(X_1)X_1]),
$$
which makes the term
$$
\frac{1}{nhf(x)}\sum_{i=1}^nK_h(x-X_i)\left(\int mf-\hat\Theta_1\right)
$$
negligible on account of condition (\ref{eq:h-negligible-3}).
\subsubsection{Proof of Lemma \ref{lem:mtg-convergence}}\label{sec:proof-mtg}
\begin{proof}[Proof of Lemma \ref{lem:mtg-convergence}] The proof is similar to \cite[Lemma 2]{WuMielniczuk2002} and \cite[Lemma 3.1]{Kulik2008EJS}.
Let $R_i=(nh\kappa_1)^{-1/2}K_h(x-X_i)\varepsilon_i/\sqrt{f(x)}$
and $\bar R_i=R_i-\Exp[R_i|{\cal F}_{i-1}]$. From the
martingale central limit theorem it suffices to show the Lindeberg
condition
$$
\sum_{i=1}^n\Exp\left[\bar R_i^21_{\{|\bar R_i|>\delta\}}\right]\to 0\quad \mbox{\rm for each } \delta>0
$$
and convergence of conditional variances
\begin{equation}\label{eq:convergence-of-cond-var}
\sum_{i=1}^n\Exp[\bar R_i^2|{\cal F}_{i-1}]\convprob 1.
\end{equation}
Let $f_{\varepsilon}$ be the density of $\varepsilon_1$. As for the
Lindeberg condition we have
\begin{eqnarray*}
\lefteqn{\sum_{i=1}^n\Exp\left[\bar R_i^21_{\{|\bar R_i|>\delta\}}\right]\le 4\sum_{i=1}^n\Exp\left[R_i^21_{\{|R_i|>\delta\}}\right]}\\
&= & C_0\frac{1}{nh}\sum_{i=1}^n\int\int K^2_h(x-u)f(u)v^2g_{\varepsilon_1}(v)1_{\{|v|>C_1\delta \sqrt{nh}\}}\\
&\le &
C_2\frac{1}{n}\sum_{i=1}^n\Exp\left[\varepsilon_i^21_{\{|\varepsilon_i|>C_1\delta\sqrt{nh}\}}\right]\to
0,
\end{eqnarray*}
where $C_0=1/(\kappa_1 f(x))$, $C_1=\left(\sqrt{C_0}\sup K(x)
\right)^{-1}$ and $C_2=C_0\int K^2$.

As for the conditional variances note first that
$$
\Exp[\bar R_i^2|{\cal F}_{i-1}]=\Exp[R_i^2|{\cal
F}_{i-1}]-\Exp\left[\left(\Exp[R_i|{\cal F}_{i-1}]
\right)^2\right]
$$
and note that the second term is of smaller order than the first
one. Now,
\begin{eqnarray*}
\lefteqn{\sum_{i=1}^n\left\{\Exp[R_i^2|{\cal F}_{i-1}]-\Exp[R_i^2]\right\}=\frac{1}{nhf(x)\kappa_1}\Exp[K_h^2(x-X_1)]\sum_{i=1}^n\left\{\Exp[\varepsilon_i^2|{\cal F}_{i-1}]-\Exp[\varepsilon_i^2] \right\}}\\
&= & \left(\frac{1}{f(x)\kappa_1}\int K^2(v)f(x-vh)dv
\right)\frac{1}{n}\sum_{i=1}^n\left\{\Exp[\varepsilon_i^2|{\cal
F}_{i-1}]-\Exp[\varepsilon_i^2] \right\}.\qquad\qquad
\end{eqnarray*}
Now, the deterministic term in the bracket is asymptotically equal
to 1. The second part converges to 0 in probability from ergodicity.
Consequently, the expression (\ref{eq:convergence-of-cond-var}) is
proven.
\end{proof}
\subsection{Proof of Propositions \ref{prop:mse_m} and \ref{lem:mse_hatm_ast}}
\begin{proof}
Recall (\ref{eq:bias}), (\ref{eq:decomposition}) and (\ref{eq:function-LRD}). Under (P1),
the result of Proposition \ref{prop:mse_m} follows from (\ref{eq:9}), (\ref{eq:10}). Under (P2), we use the expansion  with (\ref{eq:9}), (\ref{eq:7}), (\ref{eq:7a}), (\ref{eq:7b}).\\

As for Proposition \ref{lem:mse_hatm_ast}, note that all considerations for $N_{n,2}'(x)-E_n(x)$ leading to (\ref{eq:Nn2}) are in fact in $L^2$. Therefore, on account of (\ref{eq:var-assumption})
$$
\Var\left(N_{n,2}'(x)-E_n(x)\right)\sim h^4n^{-\alpha}+O\left(n^{-(\alpha+\alpha_X)}\right).
$$
The first part is of course negligible w.r.t. the bias term.

Moreover, writing
\begin{eqnarray*}
\lefteqn{\frac{1}{nhf(x)}\sum_{i=1}^nK_h(x-X_i)\left(\int mf-\hat\Theta_1\right)=\frac{1}{hf(x)}\Exp[K_h(x-X_1)]\left(\int mf-\hat\Theta_1\right)}\\
&&+\frac{1}{nhf(x)}\sum_{i=1}^n\left(K_h(x-X_i)-\Exp[K_h(x-X_1)]\right)\left(\int mf-\hat\Theta_1\right),\qquad\qquad\qquad
\end{eqnarray*}
we obtain that the variance contribution of this term is $A_1^2\Exp^2[m(X_1)X_1]n^{-\alpha_X}$.
\end{proof}
\subsection{Cross validation properties}
Under the condition (E2) one can establish the following moment bounds:
\begin{equation}\label{eq:cov-bound}
\Exp[\varepsilon_i^2\varepsilon_j\varepsilon_{j'}]=O(\Cov(\varepsilon_j,\varepsilon_{j'})),
\end{equation}
\begin{equation}\label{eq:cov-bound-1}
\Exp[\varepsilon_j\varepsilon_{j'}\varepsilon_l\varepsilon_{l'}]=O\left(\Exp[\varepsilon_j\varepsilon_{j'}]\Exp[\varepsilon_l\varepsilon_{l'}]\right),
\end{equation}
\begin{equation}\label{eq:cov-bound-2}
\Cov(\varepsilon_j^2,\varepsilon_{j'}^2)=O(\Cov(\varepsilon_j,\varepsilon_{j'})).
\end{equation}
\subsubsection{Asymptotic expansion for $\ASE(h)$}
Recall that we work under the condition (P1).
Define
\begin{equation*}
R(x)=\frac{1}{nh\hat f_h(x)}\sum_{j=1}^nK_h\left(x-X_j\right)m(X_j)-m(x).
\end{equation*}
Write
\begin{eqnarray*}
\lefteqn{\ASE(h)= I_{21}+I_{22}+I_{23}:= \frac{1}{n^3h^2}\sum_{i=1}^n\frac{1}{\hat f_h^2(X_i)}\left(\sum_{j=1}^nK_h\left(X_i-X_j\right)\varepsilon_j\right)^2}\\
&&+\frac{1}{n}\sum_{i=1}^nR^2(X_i)+\frac{2}{n}\sum_{i=1}^nR(X_i)\frac{1}{nh\hat f_h(X_i)}\sum_{j=1}^nK_h\left(X_i-X_j\right)\varepsilon_j.
\end{eqnarray*}
Let $\rho(x)=(mf)''(x)-m(x)f''(x)$. Uniformly over $\{x:f(x)>0\}$ we
have
\begin{equation}\label{eq:density-esp}
R(x)-\frac{h^2\kappa_2}{2}\frac{\rho(x)}{f(x)}=O(h^4(1+o_P(1))).
\end{equation}
Using (\ref{eq:density-esp}) we write the second part as
$$
I_{22}=\kappa_2^2\frac{h^4}{4n}\sum_{i=1}^n\left[\frac{\rho(X_i)}{f(X_i)}\right]^2 \left(1+o_P(1)\right),
$$
so that via $\Exp\left[\left(\frac{\rho(X_i)}{f(X_i)}\right)^2\right]=\int (\rho(x))^2/f(x)\; dx$ we conclude
\begin{equation}\label{eq:2}
I_{22}-\kappa_2^2\frac{h^4}{4}\int \frac{\rho^2(x)}{f(x)}\; dx=o_P(h^4/\sqrt{n}).
\end{equation}
Now, if $A_i$, $i\ge 1$, are random variables with the same mean,
and $\bar A_i=A_i-\Exp[A_i]$, then we have the following
decomposition (which will be used many times):
\begin{equation*}
\sum_{i=1}^nA_i\varepsilon_i=\Exp[A_1]\sum_{i=1}^n\varepsilon_i+\sum_{i=1}^n\bar A_i\varepsilon_i .
\end{equation*}
Typically, in LRD setting, the first part dominates. Bearing in mind
the above decomposition and since $\hat f_h$ is a consistent
estimator of $f$,
\begin{eqnarray*}
I_{23}&=&\frac{h^2}{n}\frac{1}{nh}\sum_{i=1}^n\frac{\rho(X_i)}{f^2(X_i)}\sum_{j=1}^nK_h\left(X_i-X_j\right)\varepsilon_j\left(1+o_P(1)\right)\\
&=&\Exp\left[\frac{\rho(X_1)}{f^2(X_1)}K_h\left(X_1-X_2\right)\right]\frac{h}{n}\sum_{j=1}^n\varepsilon_j\left(1+o_P(1)\right)\\
&&+\frac{h^2}{n}\frac{1}{nh}\sum_{i=1}^n\sum_{j=1}^n\overline{\frac{\rho(X_i)}{f^2(X_i)}K_h\left(X_i-X_j\right)}\varepsilon_j\left(1+o_P(1)\right).
\end{eqnarray*}
The second term is negligible w.r.t. to the first one. Noting that
\begin{equation*}
\Exp\left[\frac{\rho(X_1)}{f^2(X_1)}K_h\left(X_1-X_2\right)\right]=h\int \rho(x)\; dx+O(h^3),
\end{equation*}
we get
\begin{equation}\label{eq:3}
I_{23}=O_P\left(\frac{h^2}{n}\sum_{j=1}^n\varepsilon_j\right).
\end{equation}
It remains to deal with $I_{21}$. Split it as
\begin{eqnarray*}
\lefteqn{\frac{1}{n^3h^2}\sum_{i=1}^n\frac{1}{f^2(X_i)}\sum_{j,j'=1\atop j\not=j'}^nK_h\left(X_i-X_j\right)K_h\left(X_i-X_{j'}\right)\varepsilon_j\varepsilon_{j'}}\\
&&+\frac{1}{nh^2}\Exp\left[\frac{1}{f^2(X_1)}K_h^2\left(X_1-X_2\right)\varepsilon_1^2\right](1+o(1))\\
&&
+\frac{1}{n^3h^2}\sum_{i,j=1}^n\frac{1}{f^2(X_i)}\left(K_h^2\left(X_i-X_j\right)\varepsilon_j^2-\Exp\left[K^2_h\left(X_i-X_j\right)\varepsilon_j^2\right]\right)\\
&=:&I_{211}+I_{212}+I_{213}.
\end{eqnarray*}
Clearly,
\begin{equation}\label{eq:clearly}
I_{212}=\frac{1}{nh}\int K^2(u)\; du+ O(h/n).
\end{equation}
To deal with $I_{211}$, define
$$
T_i(h,X_j,X_{j'})=\frac{1}{h^2f^2(X_i)}K_h\left(X_i-X_j\right)K_h\left(X_i-X_{j'}\right),
$$
and note that
\begin{equation}\label{eq:integral-1}
\Exp[T_1(h,X_2,X_3)]=1+h^2\kappa_2\int f''+o(h^2).
\end{equation}
Split $I_{211}$ as
$$
\frac{\Exp[T_1(h,X_2,X_3)]}{n^3}\sum_{i=1}^n\sum_{j,j'=1\atop j\not=j'}^n\varepsilon_j\varepsilon_{j'} +\frac{1}{n^3}\sum_{i=1}^n\sum_{j,j'=1\atop j\not=j'}^n\overline{ T_i(h,X_j,X_{j'})}\varepsilon_j\varepsilon_{j'}.
$$
Variance of the last term is proportional to
\begin{equation*}
\frac{1}{n^6}\sum_{i,i'}\Cov\left(\sum_{j,j'=1\atop j\not=j'}^nT_i(h,X_j,X_{j'}),\sum_{j,j'=1\atop l\not=l'}^nT_{i'}(h,X_l,X_{l'})\right)\Exp[\varepsilon_j\varepsilon_{j'}\varepsilon_l\varepsilon_{l'}].
\end{equation*}
If all six indices are different than the term vanishes. If the indices $j,j',l,l'$ are different, but $i=i'$, then via (\ref{eq:cov-bound-1})
we obtain that the variance contribution is
$$
\frac{1}{n^6h^4}h^4nn^{2-\alpha}n^{2-\alpha}=O(1/n^{1+2\alpha}).
$$
If $j=l$ and $j'\not=l'$, together with $i=i'$ or $i\not=i'$, respectively, then via (\ref{eq:cov-bound}) the variance contribution is, respectively,
\begin{equation}\label{eq:5}
\frac{1}{n^6h^4}nh^3nn^{2-\alpha}=O\left(\frac{1}{n^{2+\alpha}h}\right),\qquad O\left(\frac{1}{n^{1+\alpha}}\right).
\end{equation}
If $j=l$ and $j'=l'$, the contribution is $O(1/(nh)^2)$. Consequently, via (\ref{eq:integral-1}) and (\ref{eq:5}), for $I_{211}$ we have
\begin{equation}\label{eq:I211}
I_{211}\sim  \frac{1+h^2\kappa_2f''}{n^2}\sum_{j,j'=1\atop j\not=j'}^n\varepsilon_j\varepsilon_{j'}+
 o_P\left(\frac{h^2}{n^2}\sum_{j,j'=1\atop j\not=j'}^n\varepsilon_j\varepsilon_{j'}\right)+O_P\left(\frac{1}{n^{(1+\alpha)/2}}\right)+O_P\left(\frac{1}{n^{1+\alpha/2}h^{1/2}}\right).
\end{equation}
Now, for $I_{213}$, its variance can be written as
$$
\frac{1}{n^6h^4}\sum_{j,j'=1}^n\Cov(\varepsilon_j^2,\varepsilon_{j'}^2)\Exp\left[\sum_{i,i'=1}^n\frac{1}{f^2(X_i)}K^2_h\left(X_i-X_j\right)\frac{1}{f^2(X_{i'})}K^2_h\left(X_{i'}-X_{j'}\right)\right].
$$
Using (\ref{eq:cov-bound-2}), one can verify that the above expression is
\begin{equation}\label{eq:6}
\frac{O(1)}{n^6h^4}(n^{2-2\alpha}\vee n)n^2h^2=O\left(\frac{1}{n^{2+2\alpha}h^2}\right)+O\left(\frac{1}{n^3h^2}\right).
\end{equation}
Via (\ref{eq:clearly}), (\ref{eq:I211}), (\ref{eq:6}),
\begin{eqnarray}\label{eq:4}
\lefteqn{I_{21}-\frac{1}{nh}\int K^2(u)\; du-\frac{1}{n^2}\sum_{j,j'=1\atop j\not=j'}^n\varepsilon_j\varepsilon_{j'}-\int f''\frac{\kappa_2h^2}{n^2}\sum_{j,j'=1\atop j\not=j'}^n\varepsilon_j\varepsilon_{j'}=O_P\left(\frac{h}{n}\right)}\nonumber\\
&&+o_P\left(\frac{h^2}{n^2}\sum_{j,j'=1\atop j\not=j'}^n\varepsilon_j\varepsilon_{j'}\right)+O_P\left(\frac{1}{n^{(1+\alpha)/2}}\right)+O_P\left(\frac{1}{n^{1+\alpha/2}h^{1/2}}\right)+O_P\left(\frac{1}{n^{1+\alpha}h}\right).\nonumber
\end{eqnarray}
Furthermore,
$$
\Exp[\frac{1}{n^2}\sum_{j,j'=1\atop j\not=j'}^n\varepsilon_j\varepsilon_{j'}]\sim C_1^2n^{-\alpha}.
$$
Combining this with (\ref{eq:2}) and (\ref{eq:3}), we have
\begin{eqnarray}\label{eq:ISE}
\lefteqn{\ASE(h)-\frac{1}{nh}\int K^2(u)\; du-\kappa_2^2\frac{h^4}{4}\int \frac{\rho^2(x)}{f(x)}\; dx-C_1^2n^{-\alpha}-C_1^2\kappa_2h^2n^{-\alpha}\int f''=}\nonumber\\
&&o_P\left(\frac{h^2}{n^2}\sum_{j,j'=1\atop j\not=j'}^n\varepsilon_j\varepsilon_{j'}\right)+O_P\left(\frac{h^2}{n}\sum_{j=1}^n\varepsilon_j\right)+O_P\left(\frac{h}{n}\right)+O_P\left(\frac{1}{n^{(1+\alpha)/2}}\right)\nonumber\\
&&+O_P\left(\frac{1}{n^{1+\alpha/2}h^{1/2}}\right)+O_P\left(\frac{1}{n^{1+\alpha}h}\right).
\end{eqnarray}
Consequently, (\ref{eq:conclusion-1}) is proven.
\subsubsection{Asymptotic expansion for $\CV(h)$}
Recall that
$$
\CV(h) :=\frac{1}{n}\sum_{i=1}^n (Y_i-\hat m_{i,h}(X_i))^2.
$$
Note that we may write
\begin{eqnarray*}
\CV(h)&=&\frac{1}{n}\sum_{i=1}^n (Y_i-\hat m_{h}(X_i))^2+O_{P}(1/(nh))\\
&=& \ASE(h)+\frac{1}{n}\sum_{i=1}^n\varepsilon_i^2+\frac{2}{n}\sum_{i=1}^n\varepsilon_i\left(m(X_i)-\hat m_h(X_i)\right)+O_{P}(1/(nh)).
\end{eqnarray*}
Here, $O_{P}(1/(nh))$ comes from replacing $\hat m_{i,h}(\cdot)$ by
$\hat m_h(\cdot)$. The second last term is treated in the very same
way as we dealt with $I_{23}$ and $I_{211}$, see (\ref{eq:3}) and
(\ref{eq:I211}), respectively. Therefore, it is
$$
O_P\left(\frac{h^2}{n}\sum_{j=1}^n\varepsilon_j+\frac{1}{n^2}\sum_{j,j'=1}^n\varepsilon_j\varepsilon_{j'}\right).
$$
Furthermore,
$
\frac{1}{n^2}\sum_{j,j'=1}^n\varepsilon_j\varepsilon_{j'}=O_P(n^{-\alpha})
$
and it dominates the term $O_P\left(\frac{1}{n^{(1+\alpha)/2}}\right)$. Consequently, via (\ref{eq:ISE}),
\begin{eqnarray*}
\frac{\CV(h)- \MISE_f(h)-\frac{1}{n}\sum_{i=1}^n\varepsilon_i^2}{ \MISE_f(h_{{\rm opt}})}&=&\frac{1}{ \MISE_f(h_{{\rm opt}})}\frac{1}{n^2}\sum_{j,j'=1}^n\varepsilon_j\varepsilon_{j'}+o_p(1),
\end{eqnarray*}
uniformly over $[B_1h_{\rm opt},B_2h_{\rm opt}]$.

\section*{Acknowledgement}
The work of the first author was supported by a NSERC (Natural
Sciences and Engineering Research Council of Canada) grant. The work
of the second author was conducted while being a Postdoctoral Fellow
at the University of Ottawa.

\bibliography{KulikLorek_RandomRegression_2010_JSPI}{}
  \bibliographystyle{plain}
 \bibliographystyle{amsalpha}
\end{document}